\documentclass[12pt,a4paper]{amsart}
\usepackage{amsmath, amssymb, bm}
\usepackage[all]{xy}

\def \sst {\scriptstyle}

\def \c {\colon}
\def \cc {\, \colon \!}

\def \q {\qquad}
\def \qq {\qquad \qquad}
\def \qqq {\qquad \qquad \qquad \qquad}
\def \bu {{\scriptscriptstyle\bullet}}
\def \skp {\medskip}
\def \ndt {\noindent}
\def \Ndt {\medskip  \noindent}
\def\shp {^{\sharp}}
\def \tilde {\raise.17ex\hbox{$\scriptstyle\mathtt{\sim}$}}   %for web address

\def \adj {\dashv}
\def \uw {\!\, \raisebox{0.3ex} {$\uparrow$}}

\def \rrw  {\; \raisebox{0.7ex}{$\longrightarrow$} \hspace{-3.7ex}
\raisebox{-0.4ex}{$\longrightarrow$} \;}

\def \lrw {\;\; \raisebox{0.4ex}{$\longleftarrow$} \hspace{-3.7ex}
 \raisebox{-0.4ex}{$\longrightarrow$} \;\;}
 
\def \lrl {\;\; \raisebox{1ex}{$\longleftarrow$} \hspace{-3.7ex}
 \raisebox{0.2ex}{$\longrightarrow$} \hspace{-3.7ex} 
 \raisebox{-0.6ex}{$\longleftarrow$} \;\;}

\def \ard {\ar@{-->}}
\def \arp {\ar@{.>}}
\def \are {\ar@{->>}}
\def \arv {\ar@{}}   %void arrow, to put a label
\def \arl {\ar@{-}}    %line
\def \arld {\ar@{--}}    %dashed line
\def \arlp {\ar@{..}}    %dotted line

\def \ti {\! \times \!}
\def \jo {{\, {\scriptstyle{\vee}} \,}}
\def \me {{\, {\scriptstyle{\wedge}}\, }}
\def \Pro {\raisebox{0.44ex}{${\mbox{\fontsize{8}{10}\selectfont\ensuremath{\prod}}}$}}
\def \Sum {\raisebox{0.44ex}{${\mbox{\fontsize{8}{10}\selectfont\ensuremath{\sum}}}$}}
\def \sub {\subset}
\def\le{\leqslant}

\def \precc {\preceq}

\def \and {\mbox{ and }}
\def \for {\mbox{ for }}

\def \The {\mbox{the }}

\def \id {{\rm id\,}}
\def \op {^{{\rm op}}}
\def \Im {{\rm Im\,}}
\def\dd {\partial}
\def\ddm {\partial^-}
\def\ddp {\partial^+}

\def \al {\alpha}

\def \th {\vartheta}
\def \si {\sigma}

\def \Set {\mathsf{Set}}
\def \Cat {\mathsf{Cat}}
\def \Top {\mathsf{Top}}
\def \pTop {{\rm p}\mathsf{Top}}
\def \dTop {{\rm d}\mathsf{Top}}
\def \cTop {{\rm c}\mathsf{Top}}
\def \d {\mathsf{d}}
\def \p {\mathsf{p}}
\def \Fl {\mathsf{F}{\rm l}\,}

\def \sing  {\{*\}}

\def \bbZ {\mathbb{Z}}
\def \bbI {\mathbb{I}}
\def \bbR {\mathbb{R}}
\def \bbS {\mathbb{S}}
\def \bbT {\mathbb{T}}
\def \uI {\uw\mathbb{I}}
\def \uR {\uw\mathbb{R}}
\def \uS {\uw\mathbb{S}}
\def \uPi {\uw\Pi}
\def \uH {\uw H}
\def \rc {{\rm c}}
\def \cI {{\rm c}\mathbb{I}}
\def \cJ {{\rm c}\mathbb{J}}
\def \cR {{\rm c}\mathbb{R}}
\def \cS {{\rm c}\mathbb{S}}

\begin{document}

\title[The topology of critical processes, I ]{The topology of critical processes, I (Processes and models)}

\author[M. Grandis]{Marco Grandis}

 \address{Marco Grandis, Dipartimento di Matematica, Universit\`a di Genova, 16146-Genova, Italy.}
 \email{grandis@dima.unige.it}



\begin{abstract}
This article belongs to a subject, Directed Algebraic Topology, whose general aim is including 
non-reversible processes in the range of topology and algebraic topology. Here, as a further 
step, we also want to cover `critical processes', indivisible and unstoppable.

	This introductory article is devoted to fixing the new framework and representing processes 
of diverse domains, with minimal mathematical prerequisites. The fundamental category and 
singular homology in the present setting will be dealt with in a sequel.
\end{abstract}

 \maketitle

\section*{Introduction}\label{Intro}

\subsection{Aims}\label{0.1}
Directed Algebraic Topology is a recent subject, dating from the 1990's. It is an extension of 
Algebraic Topology, dealing with `spaces' -- typically the {\em directed spaces} studied in 
\cite{Gr1, Gr2} -- where the paths need not be reversible, with the general aim of including 
the representation of {\em irreversible processes}.

	We want to introduce a further extension, called {\em controlled spaces}, where the 
paths need not be decomposable, in order to include {\em critical processes}, indivisible and 
unstoppable, either reversible or not.

	Taking into account transformations that cannot be stopped is an unfortunate aspect of 
our time. But there are plenty of normal events which cannot be stopped or decomposed in 
parts, like quantum effects, the onset of a nerve impulse, the combustion of fuel in a piston, the 
switch of a thermostat, the change of state in a memory cell, deleting a file in a computer, the 
action of a siphon, the eruption of a geyser, an all-or-nothing transform in cryptography, 
moving in a section of an underground network, etc.

	Critical processes and transport networks are often represented by graphs, in an effective 
way as far as they do not interact with continuous variation. We want to show that they can 
also be modelled by structured spaces, in a theory that includes classical topology and 
`non-reversible spaces'. Controlled spaces can unify aspects of continuous and discrete 
mathematics.

	It is too early to think of applications. Nevertheless, the simple fact of classifying 
phenomena of diverse domains by mathematical models which live in the same world may 
have an interest: these models can be combined together, and studied by extensions of the 
usual tools of Algebraic Topology.

	In this introductory part we fix the general framework, presenting many models and their 
concrete interpretations. The mathematical background is essentially restricted to elementary 
topology and basic category theory: limits and adjoint functors. Subsequent articles will deal 
with the homotopy theory of controlled spaces, their fundamental category and singular 
homology theory, the breaking of symmetries.

\subsection{An example}\label{0.2}
On-off controllers are systems overseeing a certain variable. Their description, as in the usual 
figure below, combines classical topology, where the variable moves freely, and graph theory, 
where a change of state takes place. We want to model them in one framework -- an enriched 
form of topology.

	For concreteness, let us think of a cooling system, with a thermostat set at temperature 
$T_0$, and a tolerance interval $[T_1, T_2]$. In the following picture the horizontal axis measures 
the temperature, and the vertical axis denotes two states: at level 0 the cooling device is off, 
at level 1 it is on

%
% FIGURE 0.2.1  On-off controller
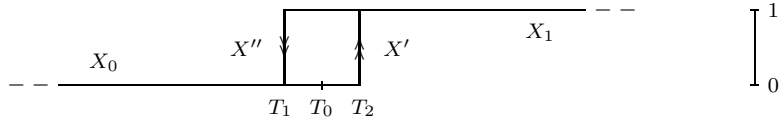
\xy <.5mm, 0mm>:
% labels + dummy pts left/above, below
(12,6) *{\sst{X_0}}; (50,10) *{\sst{X''}}; (90,10) *{\sst{X'}};  (128,14) *{\sst{X_1}}; 
(59,-6) *{\sst{T_1}}; (70,-6) *{\sst{T_0}}; (81,-6) *{\sst{T_2}}; 
(190,0) *{\sst{0}}; (190,20) *{\sst{1}};  
(-30,30) *{}; (-30,-20) *{};
% branches X_0, X_1
\POS(0,0) \arl+(80,0),  \POS(60,20) \arl+(80,0),
\POS(-13,0) \arld+(10,0),  \POS(143,20) \arld+(10,0),  
\POS(70,-1) \arl+(0,2),
% branches X', X'' and double arrows
\POS(60,0) \arl+(0,20),  \POS(80,0) \arl+(0,20),
\POS(60,12) \are+(0,-5),  \POS(80,8) \are+(0,5),
%vertical axis
\POS(185,0) \arl+(0,20), \POS(184,0) \arl+(2,0), \POS(184,20) \arl+(2,0), 
\endxy

	On the left branch $ X_0 $ the system is in stand by; if the temperature reaches $ T_2$ 
the cooling device goes on, jumping to state 1; from there, if the temperature cools to $ T_1$, 
it goes back to state 0.

	An elementary hysteresis process, or `hysteron', behaves the same way: for instance, 
the change of state in a memory cell, or the change of orientation in an elementary domain of 
a ferromagnetic object.

	We shall construct a model of this process, `pasting' two {\em natural} intervals 
$ X_0, X_1$ (with euclidean topology and nothing more) and two {\em one-jump} intervals 
$ X', X'' $ where the paths allowed have to jump the whole interval, in the marked direction: 
see \ref{3.1}(a). More complex models can be used for combined systems, like a heating and 
cooling device, in \ref{3.1}(c), or a system meant to regulate two variables, for instance 
temperature and pressure in an air-supported dome, in \ref{3.2}.

\subsection{Directed and controlled spaces}\label{0.3}
Directed spaces, our main structure meant to cover irreversible processes, was 
introduced in \cite{Gr1} and extensively studied in a book on Directed Algebraic Topology 
\cite{Gr2}; it is largely used in the theory of concurrent processes, see \ref{0.6}.

	A {\em directed space} $ X$, or {\em d-space}, is a topological space equipped with a set 
$ X\shp $ of {\em directed paths} $ [0, 1] \to X$, or {\em d-paths}, closed under: trivial loops, 
concatenation and partial increasing reparametrisation (including restrictions to subintervals). 
The selected paths, generally, cannot be travelled backwards but are {\em reflected} in the 
opposite d-space $ X\op$.

	A topological space has a {\em natural} structure of d-space, where all paths are selected. 
Directed Algebraic Topology is an extension of the classical case; in particular, the fundamental 
groupoid and the groups of singular homology are extended to directed versions: the 
fundamental category $ \uPi_1(X)$ and the preordered abelian groups $\uH_n(X) $ of 
directed singular homology.

	For all this we shall mainly refer to the book \cite{Gr2}. The prefixes {\em d-} and $\uw$ 
are used to distinguish a directed notion from the corresponding `reversible' one.

\skp	We now relax the axioms of d-spaces, to include critical processes: essentially, the 
selected paths are no longer required to be closed under restriction; they are called 
{\em controlled paths}, and the prefix {\em c-} is used to distinguish the new notions. This is still 
a directed setting, pertaining to Directed Algebraic Topology.

	In this extension we gain models of phenomena which have no place in the previous 
setting, and interesting formal `shapes', like the {\em one-jump interval} $\cI$, the 
{\em one-stop circle} $\cS^1$, the {\em n-stop circle} $\rc_n\bbS^1$, or the higher controlled 
spheres and tori described in \ref{2.3}-\ref{2.6}.

	We also loose some good properties of the theory of d-spaces. The fundamental category 
and directed singular homology of d-spaces will be extended to c-spaces, but new methods of 
computation will be needed: the van Kampen theorem and the Mayer-Vietoris sequence are 
both based on the subdivision of paths and homological chains, which is no longer permitted. 
Nevertheless, `rigid' c-spaces as the previous ones, and their fundamental category, can be 
fairly simple to analyse, precisely because of the scarcity of allowed paths.

	Essentially, the previous setting of d-spaces extends classical topology by breaking the 
symmetry of reversion: the allowed paths need no longer be reversible and the fundamental 
groupoid becomes a category. This further extension to c-spaces breaks a flexibility feature 
that d-spaces still retain: paths can no longer be subdivided, and this has drastic consequences.

\subsection{The threshold effect}\label{0.4}
As another example, in the {\em threshold effect}, or {\em siphon effect}, the process is partially 
described by a variable $v$ which can vary in a real interval $[v_0, v_1]$; when the variable 
reaches the highest value, the {\em threshold} $v_1$, it jumps down to the least value 
$ v_0$, in a way that cannot be stopped within the process itself.

 	There are many examples of this effect in Particle Physics, Natural Sciences, 
Computer Science, Medicine, Economics, Sociology, etc. 

	Some cases are well-known:

\Ndt - in Hydraulics: the empting of a basin through a siphon (see \ref{3.3});

\Ndt - in Biology: the onset of a nerve impulse ($v$ is an electric potential);

\Ndt - in Engineering: the combustion in a piston ($v$ is the quantity of fuel);

\Ndt - in Zoology and Sociology: mass migration ($v$ is the rate of the population present in a 
region, with respect to the sustainable population).

\skp	The {\em anti-siphon effect} behaves in the opposite way: the threshold is at the lowest 
level $v_0$; reaching it, the variable goes up to the highest value. The management of stocks 
of a given article, in a store or at home, roughly behaves this way.

	Two models are proposed for the siphon effect, in \ref{3.3}.

\subsection{An outline}\label{0.5}
In Section 1 we introduce the category $ \cTop $ of controlled spaces and we recall the category 
$\dTop$ of directed spaces. We also examine the links among them and other domains: 
the categories $\Top$ of topological spaces and $\pTop$ of preordered spaces. Flexible 
and rigid paths, critical paths and critical points are dealt with in \ref{1.6}.

	Section 2 begins with limits and colimits for c-spaces and d-spaces. Then we describe 
diverse c-structures on the interval $[0, 1]$, on the spheres, on the square $[0, 1]^2$, etc.

	Finally, Section 3 explores less elementary processes and how they can be modelled: 
on-off controllers in \ref{3.1} and \ref{3.2}; the threshold effect in \ref{3.3}; transport networks in 
\ref{3.4}.

\subsection{Notation and literature}\label{0.6}
The symbol $\sub$ denotes weak inclusion. A continuous mapping between topological spaces, 
possibly structured, is called a {\em map}. Open and semiopen intervals of the real line are 
always denoted by square brackets, like $]0, 1[$, $[0, 1[ $ etc. Marginal remarks are written 
in small characters.

	A {\em preorder} relation, generally written as $x \prec y$, is assumed to be reflexive and 
transitive; an {\em order} relation, often written as $x \le y$, is also assumed to be anti-symmetric. 
A mapping which preserves (resp.\ reverses) preorders is said to be {\em increasing} (resp.\ 
{\em decreasing}).

	The framework of d-spaces, its fundamental category and singular homology are used by 
various authors working in the theory of concurrency  by methods of Directed Algebraic 
Topology. This topic is covered in a recent book by L. Fajstrup, E. Goubault, E. Haucourt, 
S. Mimram and M. Raussen\cite{FjGHMR}, and many articles among which 
\cite{CaGM, FjR, Gb, GbM, MeR, Ra1, Ra2}.

	For an analysis of hysteretic processes, in the form of operators turning an input function 
into an output function, we refer to the book \cite{BrS}.

%%%
%%%
\section{Spaces with selected paths}\label{s1}

	We introduce the category $\cTop$ of controlled spaces, or c-spaces, an extension of the 
category $\dTop$ of directed spaces studied in \cite{Gr1, Gr2}, and we examine the links 
between them. Both structures are based on topological spaces with `selected paths' satisfying 
some axioms, more general for the new structure.

\subsection{Spaces and preordered spaces}\label{1.1}

	$\Top$ is the category of topological spaces and continuous mappings, or {\em maps}.

	A {\em preordered topological space} is just a space equipped with a preorder relation 
$x \prec x'$ (reflexive and transitive), without assuming any relationship between these 
structures. They form the category $\pTop$ of preordered topological spaces, with the 
increasing (i.e.\ preorder preserving) continuous mappings.

	A preordered topological space $ X $ is a `directed notion', which can be reversed: 
the object $ X\op $ has the opposite preorder $ x \prec\op x' $ (defined by $ x' \prec x$). This 
gives a (covariant) involutive endofunctor, called {\em reversor}
    \begin{equation}
R\c \pTop \to \pTop,   \qqq   RX  =  X\op.
    \label{1.1.1} \end{equation}

\begin{small}

	(The category $\Cat$ of small categories has a similar reversor.)

\end{small}

	$\bbR $ will denote the euclidean line as a topological space, and $ \bbI $ the standard 
euclidean interval $ [0, 1]$. Similarly $ \bbR^n $ and $ \bbI^n $ are euclidean spaces. 
$ \bbS^n $ is the $n$-dimensional sphere.

	On the other hand, $ \uR $ and $ \uI $ are ordered topological spaces, with their natural 
(total) order; $ \uR^n $ and $ \uI^n $ are cartesian powers in $ \pTop$, with the product order: 
$ (x_i) \le (y_i) $ if and only if, for all $ i$, $x_i \le y_i$.

	Homotopy theory in $ \Top $ is parametrised on $ \bbI$. In $\pTop $ it is parametrised on 
the ordered interval $ \uI$, yielding an elementary form of directed homotopy (cf.\ \cite{Gr2}, 
1.1.3-5).

\subsection{The terminology of paths}\label{1.2}
In a topological space $ X$, a (continuous) map $ a\c \bbI \to X $ is called a {\em path} in $ X$, 
from $ a(0) $ to $ a(1)$ -- its endpoints. It is a {\em loop} at the point $ x $ if $ a(0) = x = a(1)$.

	We begin by listing the (rather standard) terminology that we shall use for paths.

\Ndt (a) {\em Concatenation}. The concatenation of paths will be written in additive notation; the 
constant (or trivial, or degenerate) loop at the point $ x $ is written as $ 0_x$; the {\em opposite} 
path $ t \mapsto a(1 - t) $ as $ - a$.

	We recall that the (standard) concatenation $ a = a' + a'' $ of two consecutive paths 
$ a', a'' $ (with $ a'(1) = a''(0)$) is defined as
   \begin{equation} 
a(t)  =
    \begin{cases}
a'(2t),    &  \for        0 \le t \le 1/2,
\\[3pt]
a''(2t - 1),   \;  &  \for          1/2 \le t \le 1.
    \end{cases}
    \label{1.2.1} \end{equation} 

	As an important feature of topological spaces, called here the {\em path-splitting property}, 
every path $ a $ has a unique decomposition $ a = a' + a''$, with:
    \begin{equation}
a'(t)  =  a(t/2),   \qq   a''(t)  =  a((t + 1)/2))   \qq   (t \in \bbI).
    \label{1.2.2} \end{equation}

	The operation of concatenation is not associative, the constant loops do not behave as 
identities, and the opposite paths are not additive inverses -- except in trivial cases (e.g.\ in 
discrete spaces). But this works up to homotopy with fixed endpoints, which allows us to define 
the fundamental groupoid $ \Pi_1(X) $ of a space, and the fundamental group 
$ \pi_1(X, x_0) $ of a pointed space.

\Ndt (b) {\em Regular concatenation}. The {\em regular n-ary concatenation} 
$ a = a_1 + ... + a_n $ of consecutive paths is based on the regular partition 
$ 0 < 1/n < 2/n < ... < 1 $ of the standard interval, and is again uniquely determined (it is 
understood that $i = 1,..., n$):
    \begin{equation} \begin{array}{ll}
a(t)  =  a_i(nt - i + 1),   &   \for  t \in [(i - 1)/n, i/n],
\\[5pt]
a_i(t) = a((t + i - 1)/n),   \q  &   \for  t \in \bbI.
    \label{1.2.3} \end{array} \end{equation}

\Ndt (c) {\em General concatenation}. More generally, $ a = C((a_i), (t_i)) $ will denote a 
{\em general concatenation} of $ n $ consecutive paths $ a_1, ..., a_n$, based on an 
arbitrary partition $ 0 = t_0 < t_1 < ... < t_n = 1 $ of $ \bbI$
    \begin{equation} \begin{array}{lll}
a(t)  =  a_i((t - t_{i-1})/\tau_i),   &   \for  t \in [t_{i-1}, t_i],
\\[5pt]
a_i(t)  =  a(\tau_i t + t_{i-1}),   &   \for  t \in \bbI   &   (\tau_i = t_i - t_{i-1}).
    \label{1.2.4} \end{array} \end{equation}

\Ndt (d) {\em Reparametrisation}. We are interested in reparametrising the path $ a $ as 
$ a\rho\c \bbI \to X$, where the {\em reparametrisation} $ \rho\c \bbI \to \bbI $ is any increasing 
map. We speak of a {\em global reparametrisation} if $ \rho $ is surjective, that is $ \rho(0) = 0 $ 
and $ \rho(1) = 1$. We speak of an {\em invertible reparametrisation} if $ \rho $ is an 
increasing homeomorphism, or equivalently an automorphism $ \bbI \to \bbI $ of ordered sets 
(or of ordered topological spaces).

	Plainly, all $n$-ary concatenations are equivalent, up to invertible reparametrisation.

	Of course, a non-surjective reparametrisation `restricts' a path: for instance, if 
$ \rho(t) = t/2 $ (as in formula \eqref{1.2.2}), the path $ a\rho $ covers the first half of $ a$; let us 
note that it is still parametrised on $ \bbI$. More drastically, if $ \rho $ is constant $ a\rho $ is a 
constant loop.

	A {\em restriction} will be an affine, non-degenerate (i.e.\ non-constant), increasing map:
    \begin{equation}
\rho\c \bbI \to \bbI   \q   \rho(t)  =  (t_2 - t_1)t + t_1   \q   (0 \le t_1 < t_2 \le 1).
    \label{1.2.5} \end{equation}

	By the usual pleonastic terminology, a reparametrisation will also be called a 
{\em partial reparametrisation} when we want to stress that it is not assumed to be global 
(although it might be).

\subsection{Main definitions, I}\label{1.3}
A {\em controlled space} $ X$, or {\em c-space}, will be a topological space equipped with a set 
$ X\shp $ of (continuous) maps $ a\c \bbI \to X$, called {\em controlled paths}, or {\em c-paths}, 
that satisfies three axioms:

\Ndt (csp.0) ({\em constant paths}) the trivial loops at the endpoints of a controlled path are 
controlled,

\Ndt (csp.1) ({\em concatenation}) the controlled paths are closed under path concatenation: 
if the consecutive paths $ a, b $ are controlled, their concatenation $a + b$ is also,

\Ndt (csp.2) ({\em global reparametrisation}) the controlled paths are closed under 
pre-composition with every surjective increasing map $ \rho\c \bbI \to \bbI$: if $ a $ is a 
controlled path, $ a\rho $ is also.

\skp   As a consequence, the c-paths are also closed under general concatenation. The 
underlying topological space is written as $ U(X)$, or $ |X|$, and called the {\em support} 
of $X$.

	A {\em map} of c-spaces $ f\c X \to Y$, or {\em c-map}, is a continuous mapping between 
c-spaces which preserves the selected paths. Their category will be written as $ \cTop$.

	A c-space $ X $ is a directed notion. Reversing c-paths, by the involution $ r(t) = 1 - t$, 
yields the {\em opposite} c-space $ RX = X\op$, where $ a \in (X\op)\shp $ if and only if 
$ ar $ belongs to $ X\shp$. This defines the {\em reversor} endofunctor
    \begin{equation}
R\c \cTop \to \cTop,   \qq   RX  =  X\op.
    \label{1.3.1} \end{equation}

	The c-space $ X $ is {\em reversible} if $ X = X\op$, i.e.\ if its c-paths are closed under 
reversion. More generally, it is {\em reversive} if it is isomorphic to $ X\op$.

\subsection{Main definitions, II}\label{1.4}
Controlled spaces extend a structure introduced in \cite{Gr1}, also studied in 
\cite{Gr2} and elsewhere (see \ref{0.6}).

\skp	A {\em directed space} $ X$, or {\em d-space}, is equipped with a set $ X\shp $ of maps 
$ a\c \bbI \to X, $ called {\em directed paths}, or {\em d-paths}, that satisfies three axioms:

\Ndt (dsp.0) ({\em constant paths}) every trivial loop is directed,

\Ndt (dsp.1) ({\em concatenation}) if the consecutive paths $a, b$ are directed, their concatenation 
$a + b$ is also,

\Ndt (dsp.2) ({\em partial reparametrisation}) if $\rho\c \bbI \to \bbI$ is an increasing map and $a$ 
is a directed path, $a\rho$ is also.

\skp	The second axiom is the same of c-spaces (up to terminology), the others are stronger; 
every d-space is a c-space, and the notation $ X\shp $ for the set of selected paths is 
consistent. (In \cite{Gr1, Gr2} this set is written as $ dX$, a notation which has no good extension 
here.)

	A {\em map} of d-spaces, or {\em directed map}, or {\em d-map}, is a continuous mapping 
which preserves the directed paths. Their category $ \dTop $ is a full subcategory of $ \cTop$. 
The reversor endofunctor works in the same way.

\subsection{Standard intervals}\label{1.5}
The difference between these settings shows clearly in two structures of the euclidean interval 
$[0, 1]$.

\Ndt (a) In $ \dTop $ the {\em standard d-interval} $ \uI $ has for directed paths all the increasing 
maps $ \bbI \to \bbI$. It plays the role of the standard interval in the category $ \dTop $, because 
the directed paths of any d-space $ X $ coincide with the d-maps $ \uI \to X$.

	It may be viewed as the essential model of a non-reversible process, or a one-way route 
in transport networks. It will be represented as
%
% FIGURE 1.5.1, the d-interval
    \begin{equation}
\xy <.5mm, 0mm>:
% labels + dummy pt above, below
(0,-6) *{\sst{0}}; (40,-6) *{\sst{1}}; 
(0,10) *{}; (0,-20) *{};
% the interval
\POS(0,0) \arl+(40,0), \POS(0,-1) \arl+(0,2), \POS(40,-1) \arl+(0,2), 
\POS(19,0) \ar+(4,0),
\endxy
    \label{1.5.1} \end{equation}

\Ndt (b) In $\cTop $ the {\em standard c-interval} $\cI$, or {\em one-jump interval}, has the 
same support, with controlled paths the surjective increasing maps $\bbI \to \bbI$ and the 
trivial loops at 0 or 1. The controlled paths of any c-space $X$ coincide with the c-maps 
$\cI \to X$. It models a {\em non-reversible unstoppable process}, or a {\em one-way 
no-stop route} 
%
% FIGURE 1.5.2, the c-interval
    \begin{equation}
\xy <.5mm, 0mm>:
% labels + bullets + dummy pt above, below
(0,-6) *{\sst{0}}; (40,-6) *{\sst{1}}; (0,0) *{\bu}; (40,0) *{\bu}; 
(0,10) *{}; (0,-20) *{};
% the interval
\POS(0,0) \arl+(40,0), \POS(20,0) \are+(4,0)
\endxy
    \label{1.5.2} \end{equation}

\subsection{Flexible paths and critical points}\label{1.6}
(a) {\em Flexible paths}. In a c-space $ X$, a point $ x $ will be said to be {\em flexible} 
if its trivial loop $ 0_x $ is controlled; the {\em flexible support} $ |X|_0 $ is the subspace of 
these points. In a diagram, an isolated flexible point will be marked by a bullet, as in  
figure \eqref{1.5.2} above.

	We say that a controlled path $ a $ is {\em splittable} if its halves $ a', a'' $ (cf.\ \eqref{1.2.2}) 
are also controlled, so that the decomposition $ a = a' + a'' $ stays within c-paths; we say that 
$ a $ is {\em flexible} if all its restrictions are controlled (see \eqref{1.2.5}), or equivalently all its 
decompositions in general concatenations give raise to c-paths. Each controlled trivial loop is 
flexible. A c-map preserves all these properties.

	The c-space itself is {\em flexible} if every point and every c-path is flexible. A c-space is a 
d-space if and only if it is flexible, if and only if every trivial loop is controlled and all its 
controlled paths are splittable.

	A c-path $ a $ is {\em rigid} if in each general concatenation of $a$ by controlled paths, 
precisely one of them is not constant. A c-space is {\em rigid} if every non-trivial path is a general 
concatenation of rigid paths. The interval $ \cI $ is rigid, as well as many c-spaces introduced in 
the next section.

\Ndt (b) {\em Critical paths and critical points}. In a c-space $ X$, a controlled path is {\em critical} 
if it is not flexible.

	A point $ x $ is:

\Ndt - {\em critical}, if every non-trivial c-path $ a $ through $ x $ (i.e.\ $ x \in \Im a$) is critical, 
and there is some,

\Ndt - {\em future critical}, if every non-trivial c-path starting there is critical, and there is some,

\Ndt - {\em past critical}, if every non-trivial c-path arriving there is critical, and there is some.

\skp	A future or past critical point $ x $ is always flexible, a critical point need not. A d-space 
has no critical points.

	In the interval $ \cI $ all points are critical, the point 0 is also future critical, while 1 is also 
past critical. There are c-spaces where these three kinds are disjoint: see \ref{2.3}(e).

\subsection{Reshaping and generated structures}\label{1.7}
The c-structures on a topological space $ X $ are closed under arbitrary intersection (as subsets 
of $ \Top(\bbI, X))$, and form a complete lattice for the inclusion: we say that the structure 
$ X_1 $ is {\em finer} than $ X_2 $ if $ X_1\shp \sub X_2\shp$, or equivalently if the identity map 
of $ X $ gives a map $ X_1 \to X_2$; this map is called a {\em reshaping}.

\Ndt (a) Every set $ S $ of paths in the space $X$ generates a c-structure, the finest, or smallest, 
containing it. It is obtained by adding all the constant loops at the endpoints of the paths of $S$, 
and stabilising the latter under global reparametrisation and general concatenation.

\Ndt (b) Similarly, the d-structures on a topological space $X$ form a complete lattice. Every 
set of paths of $ X $ generates a d-structure.

\Ndt (c) If we start from a c-space $X$, the d-structure generated by the c-paths can 
be obtained stabilising them under constant paths, restriction and general concatenation.

\Ndt (d) The forgetful functor $U\c \dTop \to \Top$ takes a d-space to its support, the underlying 
topological space $|X|$. It has a left and a right adjoint
    \begin{equation}
U\c  \dTop   \lrl   \Top  \cc D, D'   \qq   D \adj U \adj D'. 
    \label{1.7.1} \end{equation}

	For a topological space $T$, the d-space $DT$ is the same space with the {\em discrete 
d-structure} (the finest, or smallest), with directed paths all the trivial loops. $D'T$ has the 
{\em indiscrete d-structure} (the largest, or coarsest), where all paths are controlled.

\Ndt  (e) The category $\cTop$ has two forgetful functors to topological spaces
    \begin{equation}
U\c \cTop \to \Top,   \qq   U_0\c \cTop \to \Top,
    \label{1.7.2} \end{equation}
where $ U(X) = |X| $ is the topological support and $ U_0(X) = |X|_0 $ is the flexible support. 
$ U $ has both adjoints, $ U_0 $ has only the left one
    \begin{equation}
D_c \adj U \adj D',   \qq\q   D \adj U_0.
    \label{1.7.3} \end{equation}

	For a topological space $T$, the c-space $D_cT$ is the same space with the 
{\em discrete c-structure}: no path is controlled. $D'T$ has the {\em indiscrete 
c-structure}, where all paths are controlled. In $DT$ all trivial loops are controlled. The 
functors $D$ and $D'$ take values in $\dTop$, and are denoted as previously.

	A topological space will be viewed as a c-space (and a d-space) by its {\em natural} 
structure $ D'T$, so that all its paths are selected.

\Ndt  (f) The singleton has two structures in $ \cTop$: the {\em c-discrete} singleton 
$ D_c\sing $ and the {\em flexible singleton} $ \sing_0 = D\sing = D'\sing, $ with a controlled 
loop $ 0_*$; the flexible singleton is by far more important, as it is the terminal object and the 
unit of the cartesian product (see \ref{2.1}).

	A d-map $x\c D_c\sing \to X$ is `the same' as a point of $X$, while a d-map 
$x\c \sing_0 \to X$ is a flexible point. In other words, $D_c\sing$ represents the functor 
$U\c \cTop \to \Set$, while $\sing_0$ represents $U_0\c \cTop \to \Set$.

\Ndt  (g) All the c-spaces $DT$ are trivially flexible and rigid.

\subsection{Comparing directed structures}\label{1.8}
We are considering three ways of enriching topological spaces by a directed structure (in a 
general sense), encoded in the categories $ \pTop$, $\dTop $ and $ \cTop$. We now examine 
their interplay.

	A preordered topological space $ X $ (in the sense recalled in \ref{1.1}) will always be 
viewed as a d-space (and a c-space) by selecting the increasing (continuous) paths 
$ \uI \to X$, where $ \uI $ denotes the ordered euclidean interval $ [0, 1]$.

	This defines a functor $ \d\c \pTop  \to \dTop$, and our categories are linked by three 
obvious functors
    \begin{equation}
\d\c \pTop \to \dTop,   \q   \dTop  \sub  \cTop,   \q   \d\c \pTop \to \cTop.
    \label{1.8.1} \end{equation}

	Let us note that $ \d $ {\em is not an embedding}: trivially, all preorders on a discrete 
topological space give the same selected paths, namely the trivial loops. (One can find 
more interesting examples in \cite{Gr2}, 1.4.5.)

\Ndt (a) There is an adjunction
    \begin{equation}
\d\c  \pTop  \lrw   \dTop  \cc \p,   \qq   \p \adj \d,
    \label{1.8.2} \end{equation}
where the left adjoint $ \p $ provides a d-space with the {\em path-preorder} $ x \precc x'$, 
meaning that there exists a d-path from $ x $ to $ x'$. The counit on a preordered space $ X $ 
is the preorder-reshaping $\p\d X \to X$: if $ x \precc x' $ there exists a d-path from $ x $ to $ x' $ 
in $\d X$, whence $ x \prec x' $ in $ X$.

	Both functors $ \p, \d $ are faithful. A d-space is said to be {\em of (pre)order type} if it can 
be obtained, as above, from a topological space with such a structure. Thus $ \uR^n $ and 
$ \uI^n $ are of order type; $ \bbR^n$, $ \bbI^n $ and $ \bbS^n $ are of chaotic-preorder type. 
The directed sphere $ \uS^n $ described in \ref{2.5} is not of preorder type (for $ n > 0$). 

\Ndt (b) The embedding $ \dTop \to \cTop $ has a left and a right adjoint:
    \begin{equation} \begin{array}{ccr}
\hat{\;}\c  \cTop \to \dTop   &\qq&   (\The reflector),
\\[5pt]
\Fl  \c  \cTop \to \dTop   &&   (\The coreflector).
    \label{1.8.3} \end{array} \end{equation}

	For a c-space $ X$, the {\em generated d-space} $\hat{X}$ has the same underlying 
topological space with the d-structure generated by the c-paths. The unit of the adjunction is 
the reshaping $ X \to \hat{X}$, the counit is the identity $\hat{Y} = Y$ for a d-space $ Y$.

	In the second construction the {\em flexible part} $\Fl X$ is the flexible support 
$|X|_0$ with the d-structure of the flexible c-paths. The counit is the inclusion $\Fl X \to X$, the 
unit is the identity $ Y = \Fl Y$ for a d-space $ Y$.

\begin{small}

\skp	The full subcategory of reversible c-spaces has a similar reflector and coreflector: the 
{\em generated reversible c-space} and the {\em reversible part}.

\end{small}

\Ndt (c) Composing the adjunction \eqref{1.8.2} with the previous reflection
    \begin{equation} \begin{array}{c} 
    \xymatrix  @C=20pt @R=20pt
{
~\pTop~  \ar@<-2pt>[r]_-{\d}    &   ~\dTop~    \ar@<-2pt>[l]_-{\p}  \ar@<-2pt>[r]_-{\sub}  &
~\cTop~     \ar@<-2pt>[l]_-{\hat{\;}}  \;&\;
~\pTop~  \ar@<-2pt>[r]_-{\d}    &   ~\cTop~     \ar@<-2pt>[l]_-{\hat{\p}}
}
    \label{1.8.4} \end{array} \end{equation}
we get the adjunction $\hat{\p} \adj \d$, where $ \d $ still equips a preordered space 
$ X $ with the increasing maps $ \uI \to X $ as c-paths (producing a d-space), while 
$\hat{\p}(X) = \p(\hat{X})$ provides a c-space with the {\em generated-path preorder} 
$ x \precc x'$, depending on the d-paths of $\hat{X}$. (The c-paths of $ X $ give a preorder 
on the flexible support $ |X|_0$, not used here.)

%%%
%%%
\section{Limits, colimits and structural models}\label{s2}
	Limits and colimits, for c-spaces and d-spaces, are easily obtained as topological limits 
and colimits with the initial or terminal structure determined by the structural maps. 

	Then we describe diverse c-structures on the interval, the spheres and the square; they 
can represent elementary events and will be used as bricks to form models of more complex 
processes.

\subsection{Limits and colimits}\label{2.1}
We already remarked that the c-structures on a topological space $ T $ form a complete lattice. 
Therefore every family of maps $f_i\c T \to X_i$ with values in c-spaces defines an initial 
c-structure on the space $T$: a path $a$ is controlled if and only if all composites $f_ia$ are.
Dually, every family of maps $ f_i\c X_i \to T $ defined on c-spaces gives raise to a final 
c-structure on the space $ T$: the controlled paths in $ T $ are generated by all the paths 
$ f_i a$, where $ a \in X_i\shp $ for some index $ i$.

	A (controlled) {\em subspace} $ X' \sub X $ of a c-space $ X $ has the initial structure of 
the embedding, which selects those paths in $ X' $ that are controlled in $ X$. A (controlled) 
{\em quotient} $ X/R $ has the quotient structure, that is the final one for the projection 
$ p\c X \to X/R$; it is generated by the projected c-paths through general concatenation (see 
\ref{1.7}(a)).

	The category $ \cTop $ has all limits and colimits, constructed as in $ \Top $ and 
equipped with the initial or final c-structure for the structural maps. For instance a path 
$\bbI \to \Pro X_i $ with values in a product of c-spaces is controlled if and only if all its 
components $ \bbI \to X_i $ are, while a path $ \bbI \to \Sum X_i $ with values in a sum is 
controlled if and only if it is in some summand $ X_i$. Equalisers and coequalisers are 
realised as subspaces or quotients, in the sense described above.

	We already described the terminal $ \sing_0$, which is the unit of the cartesian product. 
On the other hand, $ X \ti D_c\sing $ is the discrete c-structure $ D_c|X|$ on the underlying 
space.

	If $ X $ is a c-space and $ A \sub |X| $ is a {\em non-empty} subset, $ X/A $ will denote the 
c-quotient of $ X $ which identifies all points of $ A$.

\skp	All this works in the same way in $ \pTop $ and $ \dTop$. The embedding 
$ \dTop \sub \cTop $ preserves all limits and colimits, as it has both adjoints (see \eqref{1.8.3}). 
On the other hand, the canonical functors $ \d\c \pTop \to \dTop $ and $ \d\c \pTop \to \cTop $ 
of \eqref{1.8.1} preserve limits (as right adjoints) and sums (obviously), but do not preserve 
coequalisers.

	In fact, in $ \pTop $ the coequaliser of the endpoints $ \sing \rrw \uI $ is the circle 
$ \bbS^1 $ with the indiscrete preorder. In $ \dTop $ (and $ \cTop$) we get a non-trivial 
d-structure, the directed circle $ \uS^1$, described below in \ref{2.5}(a). Essentially, this is 
`why' directed homotopy is simple but very elementary in $ \pTop$.

	(The standard c-circle $\cS^1$, described in \ref{2.6}(a), is the coequaliser of the 
endpoints in $\cI$.)

\subsection{Controlled actions}\label{2.2}
	Let $G$ be a group, written in additive notation (commutative or not). A {\it controlled 
$G$-space} is a c-space $X$ equipped with a (right) {\em action} of $ G$: this is an action on 
the underlying topological space such that, for each $ g \in G$, the induced map
    \begin{equation}
X \to X,   \qq   x  \mapsto  x + g,
    \label{2.2.1} \end{equation}
is a map of c-spaces (and therefore an isomorphism of $ \cTop$). Directed $G$-spaces are a 
particular case.

	The c-space {\em of orbits} $ X/G $ is the quotient c-space, modulo 
the equivalence relation which collapses each orbit to a point. Its c-paths are simply the 
projections of the directed paths of $ X$, as verified below. The same holds for d-spaces.

\begin{small}

\skp	We have to prove that these projections are closed under global (resp.\ partial) 
reparametrisation and binary concatenation. The first fact is obvious. As to the second, let 
$ a, b\c \bbI \to X $ be two controlled paths whose projections are 
consecutive in $ X/G$: there is some $ g \in G $ such that $a(1) = b(0) + g$. Then the path 
$b'(t) = b(t) + g$ is controlled in $ X$, and $a + b'$ is also. Finally, writing $p\c X \to X/G$ 
the canonical projection, $ pa + pb = p(a + b')$ is the projection of a controlled path.

\end{small}

\subsection{Elementary models}\label{2.3}
(a) The euclidean interval $ \bbI $ and the euclidean line $ \bbR $ have the natural d-structure, 
where all paths are selected. The same holds for their cartesian powers $ \bbI^n $ and 
$ \bbR^n$, and for all spheres $ \bbS^n$. $\bbI $ will be called the {\em natural} interval.

\Ndt (b) The ordered euclidean interval $ \uI $ and the ordered euclidean line $ \uR $ have the 
d-structure given by the increasing paths (already recalled for the former). The same holds for 
their cartesian powers $ \uI^n $ and $ \uR^n$. They are not reversible (for $ n > 0$), yet 
reversive, i.e.\ isomorphic to the opposite.

	$\uI $ is the standard ordered interval, and also the standard d-interval, as already said. 
Its d-structure is generated by the identity map $ \bbI \to \bbI$; a d-map $ \uI \to X $ is the same 
as a directed path of $ X$.

	But $ \uI $ will also be important in $ \cTop$, as the {\em flexible interval}. Indeed, for a 
c-space $ X$, the c-maps $ \uI \to X $ are the flexible paths of $ X$.

\Ndt (c) We already introduced the standard controlled interval $ \cI$, with the c-structure 
generated by the identity map $ \bbI \to \bbI$: the c-paths are the surjective increasing maps 
$ \bbI \to \bbI $ and the trivial loops at the endpoints. The c-maps $ \cI \to X $ are the 
selected paths of the c-space $ X $ (possibly a d-space).

	The c-space $ \cI $ will also be called the {\em quantum interval}, or the 
{\em one-jump interval}
%
% FIGURE 2.3.1 = figure 1.5.2, the c-interval
    \begin{equation}
\xy <.5mm, 0mm>:
% labels + bullets + dummy pt above, below
(0,-6) *{\sst{0}}; (40,-6) *{\sst{1}}; (0,0) *{\bu}; (40,0) *{\bu}; 
(0,10) *{}; (0,-20) *{};
% the interval
\POS(0,0) \arl+(40,0), \POS(0,-1) \arl+(0,2), \POS(40,-1) \arl+(0,2), 
\POS(20,0) \are+(4,0)
\endxy
    \label{2.3.1} \end{equation}

	The generated d-space is $ (\cI)\hat{\;} = \uI$, while the flexible part $\Fl(\cI) = D\{0, 1\}$ is 
the discrete boundary $ \dd \bbI $ of the interval, with its trivial loops.
 
\Ndt (d) The {\em line with integral stops} $ \cR$, or {\em integral jumps}, is equipped with the 
c-structure generated by the family of embeddings $ \bbI \to \bbR$, $t \mapsto t + k $ 
($k \in \bbZ$). Now the c-paths are the increasing maps $ \bbI \to \bbR $ whose image is 
precisely an interval $ [k, k'] $ with integral endpoints (possibly the same)
%
% FIGURE 2.3.2, the c-line
    \begin{equation}
\xy <.5mm, 0mm>:
% labels + bullets + dummy pt above, below
(0,-6) *{\sst{-1}}; (40,-6) *{\sst{0}}; (80,-6) *{\sst{1}}; (120,-6) *{\sst{2}}; (160,-6) *{\sst{3}}; 
@i@={(0,0), (40,0), (80,0), (120,0), (160,0)} @@{*{\bu}};
(0,15) *{}; (0,-20) *{};
% the line
\POS(0,0) \arl+(160,0), \POS(-13,0) \arld+(10,0), \POS(163,0) \arld+(10,0), 
\POS(20,0) \are+(4,0), \POS(60,0) \are+(4,0), \POS(100,0) \are+(4,0), \POS(140,0) \are+(4,0),
\endxy
    \label{2.3.2} \end{equation}

	The line $\cR$ is a controlled $\bbZ$-space, with respect to the action of the group $\bbZ$ 
by translations. The interval $\cI$ is a subspace of $\cR$, and the latter is the 
controlled $\bbZ$-space generated by the embedding of $\cI$.

	The line $\cR$ is a rigid c-space (see \ref{1.6}): the rigid paths are those of length 1, and 
every non-trivial c-path is a concatenation of them, on a suitable partition. All points of $ \cR $ 
are critical; the integral numbers are also past and future critical. The generated d-space 
$ (\cR)\hat{\;} = \uR $ is of order type; the flexible part $\Fl (\cI) = D\bbZ$ is the discrete integral 
line with the discrete d-structure.

\Ndt(e) Let $ X $ be the euclidean ordered interval $ [0, 3]$, with controlled paths given by the 
increasing maps $ \bbI \to X $ whose image either contains the open subinterval $ ]1, 2[ $ or 
does not meet it.

	Loosely speaking, we are modelling a process measured on the interval $[0, 3]$, which

\Ndt - can only proceed `forward',

\smallskip \ndt - passing point 1, is obliged to go on to point 2, at least,

\Ndt or a one-way route with a no-stop section, or a stream with rapids
%
% FIGURE 2.3.3, a one-way route with a no-stop section
    \begin{equation}
\xy <.5mm, 0mm>:
% labels + dummy pt above, below
(0,-6) *{\sst{0}}; (40,-6) *{\sst{1}}; (80,-6) *{\sst{2}}; (120,-6) *{\sst{3}};
(0,15) *{}; (0,-25) *{};
% the line
\POS(0,0) \arl+(120,0),
\POS(0,-1) \arl+(0,2), \POS(40,-1) \arl+(0,2), 
\POS(80,-1) \arl+(0,2), \POS(120,-1) \arl+(0,2),
\POS(19,0) \ar+(4,0), \POS(60,0) \are+(4,0), \POS(99,0) \ar+(4,0),
\endxy
    \label{2.3.3} \end{equation}

	The point 1 is future critical; all the points of $ ]1, 2[$ are not flexible and critical; 2 is 
past critical. The generated d-space $\hat{X} = \uw [0, 3]$ is the ordered structure, while 
$\Fl X$ is the ordered structure on the flexible support $ [0, 1] \cup [2, 3]$.

\subsection{Other structures on the interval. }\label{2.4}
We have already seen three c-structures on the euclidean interval $ [0, 1]$: the natural structure 
$ \bbI$, where all paths are controlled; the ordered structure $ \uI$, with the increasing paths; 
the one-jump structure $ \cI$, with the surjective increasing paths and the trivial loops at the 
endpoints.

	There are many others, that can be used as bricks of modelisation. We list here some of 
them; two `siphon structures' can be found in \ref{3.3}.

\Ndt (a) The {\em two-jump interval} $ \cJ $ has the c-structure generated by the restrictions to 
the first or second half
    \begin{equation}
c^-(t)  =  t/2,   \qq   c^+(t)  =  (t + 1)/2   \qqq   (t \in \bbI),
    \label{2.4.1} \end{equation}
%
% FIGURE 2.4.1, the two-jump interval  \cJ
\xy <.5mm, 0mm>:
% labels + bullets + dummy pt above, below
(0,-6) *{\sst{0}}; (40,-6) *{\sst{1/2}}; (80,-6) *{\sst{1}};  
(0,0) *{\bu}; (40,0) *{\bu}; (80,0) *{\bu};  
(-50,10) *{}; (0,-20) *{};
% the line
\POS(0,0) \arl+(80,0), 
\POS(20,0) \are+(4,0), \POS(60,0) \are+(4,0),
\endxy

	The non-trivial c-paths are the increasing maps $ \bbI \to \bbI $ whose image is either 
$ [0, 1/2]$, or $ [1/2, 1]$, or $ [0, 1]$. This c-space is isomorphic to the subspace 
$ \rc[0, 2] \sub \cR$, and can model a {\em non-reversible two-stage process}. Formally, 
it parametrises the concatenation of two c-paths (see Part II).

\Ndt (b) {\em Delayed intervals}. The {\em past-delayed c-interval} $\rc_-\bbI$ will be the 
euclidean interval $ [0, 1] $ with the c-structure generated by the past-delayed 
reparametrisation $ \rho\c \bbI \to \bbI $
    \begin{equation}
\rho(t)  =  0 \jo (2t - 1),   \qq   \si(t)  =  2t \me 1,
    \label{2.4.2} \end{equation}
while the {\em future-delayed c-interval} $ \rc_+\bbI $ has the c-structure generated by the 
future-delayed reparametrisation $ \si\c \bbI \to \bbI$.

	In $ \rc_-\bbI $ the non-trivial controlled paths are the surjective increasing maps 
$ \bbI \to \bbI $ which are constant on some non-degenerate interval $ [0, t_1]$. For a c-space 
$ X$, a c-map $ \rc_-\bbI \to X $ is the same as a {\em past-delayed c-path} (constant as above).
	
	These c-spaces are not reversive, but anti-isomorphic to each other, by reversion
    \begin{equation}
r\c \rc_-\bbI \to (\rc_+\bbI)\op,   \qq   r(t)  =  1 - t.
    \label{2.4.3} \end{equation}

	These structures are rigid and finer than $\cI$, because their generators are 
surjective increasing maps, i.e.\ controlled paths in $\cI$.
	
	There are many delayed structures on the interval. They can model
 {\em irreversible non-stoppable processes with inertia}, or {\em inductance}

\Ndt (c) The {\em reversible d-interval} $\, \bbI^\sim \,$ is the d-space generated by identity 
and reversion $ \id, r\c \bbI \to \bbI$. Its directed paths are the piecewise monotone maps 
$ \bbI \to \bbI $ (increasing or decreasing on each subinterval of a finite partition of its domain). 
Its interest comes out of the fact that, for a d-space $ X$, the 
d-maps $ a\c \bbI^\sim \to X $ are precisely the reversible directed paths of $ X $ 
(\cite{Gr2}, 1.4.6). This d-spaces is strictly finer than the natural interval $ \bbI$. It can model 
a {\em shock absorber}.

\Ndt  (d) The {\em one-jump reversible interval} $ \cI^\sim$ has the c-structure generated 
by identity and reversion $ \id, r\c \bbI \to \bbI$. Every c-path has an integral length and those 
of length 1 are rigid; our c-space is rigid and the generated d-space is $ \bbI^\sim$. 

	The interval $ \cI^\sim $ can (basically) model a {\em reversible non-stoppable 
process}, like the change of state in a memory cell, or a {\em two-way no-stop route}, or the 
flights of a given airplane between two airports
%
% FIGURE 2.4.5 the interval   \cI^\sim
    \begin{equation}
\xy <.5mm, 0mm>:
% labels + bullets + dummy pt above, below
(0,-6) *{\sst{0}}; (50,-6) *{\sst{1}}; (80,-2) *{\cI^\sim}; (0,0) *{\bu}; (50,0) *{\bu}; 
(0,10) *{}; (0,-20) *{};
% the interval
\POS(0,0) \arl+(50,0), \POS(0,-1) \arl+(0,2), \POS(50,-1) \arl+(0,2), 
\POS(18,0) \are+(-4,0) \POS(34,0) \are+(4,0)
\endxy
    \label{2.4.5} \end{equation}

\subsection{Directed spheres and tori}\label{2.5}
(a) The {\em standard d-circle} $ \uS^1 $ is the standard circle with the {\em anticlockwise 
structure}, where the directed paths $ a\c \bbI \to \bbS^1 $ move this way, in the oriented plane 
$ \bbR^2$: $a(t) = (\cos \th (t), \sin \th (t))$, with an increasing (continuous) argument 
$ \th \c \bbI \to \bbR $
%
%  FIG. 2.5.1 (from DAT 14.3.2): the directed circle
    \begin{equation} 
\xy <1pt, 0pt>:
%  label and dummy pt above, below
(45,-12) *{\uS^1}; (0,30) *{}; (0,-30) *{};
%
% circle & arrow
(0,0) = "A"  *\cir<25pt>{}="ca"  
\POS(1,23.8)\ar-(2,0)
\endxy
    \label{2.5.1} \end{equation}

	The directed circle can be described as an orbit space
    \begin{equation}
\uS^1  =  \uR/\bbZ,
    \label{2.5.2} \end{equation}
with respect to the action of the group of integers on the directed line $\uR$, by translations: 
the directed paths of $\uS^1$ are simply the projections of the increasing paths in the line.

	The c-space  $\uS^1$ can also be obtained as the coequaliser in $ \dTop $ of the 
following pair of maps
    \begin{equation}
\ddm, \ddp\c \sing  \rrw   \uI,   \qq   \ddm(*)  =  0,   \q   \ddp(*)  =  1.
    \label{2.5.3} \end{equation}

\skp	Indeed, this coequaliser is the quotient $ \uI/\dd \bbI$, which identifies the endpoints; 
the d-structure of the quotient, generated by the projected paths, is what we want: it is sufficient 
to concatenate a finite number of projected paths (which are already stable under partial 
reparametrisation).

\Ndt (b) The {\em standard directed n-sphere} is defined, for $ n > 0$, as the quotient of the 
directed cube $ \uI^n $ modulo its (ordinary) boundary $ \dd \bbI^n$
    \begin{equation}
\uS^n  =  (\uI^n)/(\dd\bbI^n),   \q\;\;\;   \uS^0  =  \bbS^0  =  \{-1, 1\}   \qq   (n > 0).
    \label{2.5.4} \end{equation}
while $ \uS^0 $ has the discrete topology and the natural d-structure (obviously discrete)

	All directed spheres are reversive; their d-structure can be described by an asymmetric 
distance (see \cite{Gr2}, 6.1.5). The pointed suspension of $ \bbS^0 $ in the category of 
pointed d-spaces gives $ \uS^1 $ and, by iteration, all higher $ \uS^n $ (\cite{Gr2}, 1.7.4, 1.7.5). 
The unpointed suspension gives different d-spaces, less interesting.

\Ndt (c) The {\em standard directed n-torus} is a cartesian power of $\uS^1$ 
    \begin{equation}
\uw\bbT^n  =  (\uS^1)^n.
    \label{2.5.5} \end{equation}

	Equivalently, it is the orbit d-space $(\uR^n)/\bbZ^n$, for the action of 
the additive group $\bbZ^n$ by translations.

\subsection{Controlled spheres and tori}\label{2.6}
(a) The {\em standard c-circle} $ \cS^1$, or {\em one-stop circle}, is now defined as the 
following orbit c-space
%
%  FIG. 2.6.1: the standard controlled circle
    \begin{equation} 
\xy <1pt, 0pt>:
%  bullet, label and dummy pt above, below
(23.8,0) *{\bu}; (80,-12) *{\cS^1  =  (\cR)/\bbZ}; (0,30) *{}; (0,-30) *{};
%
% circle & double arrow
(0,0) = "A"  *\cir<25pt>{}="ca"  
\POS(0,23.8)\are-(4,0)
\endxy
    \label{2.6.1} \end{equation}
for the action of the group $\bbZ$ on the line $ \cR$, by translations. The controlled paths of 
$ \cS^1 $ are simply the projections of the controlled paths in the line: here this means an 
anticlockwise path (as in \ref{2.5}(a)) which is a loop at $ [0]$, the only flexible point 
(corresponding to $ (1, 0) $ in the plane). The simple loops are rigid, and so is $ \cS^1$.

	The circle $\cS^1 $ can also be obtained as the coequaliser in $ \cTop $ of the endpoints 
of $ \cI$
    \begin{equation} 
	\ddm, \ddp\c \sing   \rrw   \; \cI,   \qq   \ddm(*)  =  0,  \q  \ddp(*)  =  1.
    \label{2.6.2} \end{equation}

	All points are critical; the flexible point is also past and future critical. The generated 
d-space is $ (\cS^1)\hat{\;} = \uS^1$, while $\Fl (\cS^1)$ is the base-point with its trivial loop.

\Ndt(b) More generally, the {\em n-stop c-circle} $\rc_n\bbS^1 $ ($n > 0$) is the orbit space
    \begin{equation} 
\rc_n\bbS^1  =  (\rc_n\bbR)/\bbZ   \qq   (\rc_1\bbS^1  =  \cS^1),
    \label{2.6.3} \end{equation}
where the c-paths of $\rc_n\bbR $ are the increasing paths whose image is an interval 
$ [k/n, k'/n]$, for integers $ k \le k'$. 

	In $ \rc_n\bbS^1 $ a c-path is an anticlockwise path between two points $ [k/n] $ and 
$ [k'/n] $ of the circle. The `minimal generators' have length $ 1/n $ of the circle and are rigid; 
the c-space itself is also.

\begin{small}

\skp	Rotating motions can follow this pattern, with mandatory direction and stops: for instance, 
the second hand of a watch, a washing machine dial, a vertical panoramic wheel with $ n $ 
cabins. The mode dial of a photocamera and a railway turntable can be modelled by the 
reversible c-space generated by $ \rc_n\bbS^1$.

\end{small}

\Ndt (c) The {\em standard c-sphere} $ \cS^n $ is defined as a quotient of the cube $ \cI^n $ 
(for $ n > 0$)
    \begin{equation} 
\cS^n  =  (\cI^n)/(\dd\bbI^n),   \q\;\;\;   \cS^0  =  \bbS^0  =  \{-1, 1\}   \qq   (n > 0).
    \label{2.6.4} \end{equation}

	It will be examined in Part II, were we shall see that the pointed suspension of 
$ \bbS^0 $ in the category of pointed c-spaces gives all $ \cS^n$. (Again, the unpointed 
suspension gives different c-spaces.)

\Ndt (d) The {\em standard controlled n-torus} is a cartesian power of $\cS^1$ 
    \begin{equation}
\rc\bbT^n  =  (\cS^1)^n,
    \label{2.6.5} \end{equation}
and can also be obtained as the orbit c-space $(\cR^n)/\bbZ^n$.

\subsection{Controlled squares and cubes}\label{2.7}
\Ndt (a) We have already seen the ordered square $ \uI^2$, also called the {\em flexible square} 
when viewed in $ \cTop$.

\Ndt (b) The {\em standard c-square} $ \cI^2$, represented in the left figure below, has the 
structure of a cartesian power: a path $ \bbI \to \bbI^2 $ is controlled if and only if it is 
increasing and each of its projections covers $ [0, 1]$, or is constant at  0  or  1
%
% FIGURE 2.7.1, Standard controlled square and the generated d-space
    
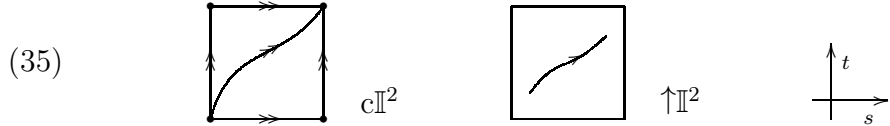
\begin{equation} 
\xy <.5mm, 0mm>:
%horizontal POSITIONS of two squares: [0,30], [80,110]
%
% dummy pt above/left, below
(-20,22) *{}; (0,-22) *{};
% labels, bullets
(45,-10) *{\cI^2}; (125,-10) *{\uI^2}; (175,-15) *{\sst{s}}; (169,0) *{\sst{t}};
(0,-15) *{\bu}; (30,-15) *{\bu}; (30,15) *{\bu};  (0,15) *{\bu}; 
% two squares
@i@={(0,-15), (30,-15), (30,15), (0,15)},
s0="prev"  @@{;"prev";**@{-}="prev"};
@i@={(80,-15), (110,-15), (110,15), (80,15)},
s0="prev"  @@{;"prev";**@{-}="prev"};
% curved paths and their arrows
(0, -15); (30, 15) **\crv{(5,5)&(20,0)},
(85, -8); (105, 7) **\crv{(92,2)&(95,-3)},
\POS(15,3) \are+(4,2), \POS(95,0) \ar+(4,2), 
% double arrows
\POS(15,-15) \are+(4,0), \POS(15,15) \are+(4,0),  
\POS(0,0) \are+(0,4), \POS(30,0) \are+(0,4),
% axes of the plane
\POS(160,-10) \ar+(20,0),  \POS(165,-15) \ar+(0,20),
\endxy
    \label{2.7.1} \end{equation}

	There are four flexible points, the vertices of the square. The c-paths of $ \cI^2 $ have five 
kinds of generators: two horizontal paths $ s \mapsto (s, \al) $ (for $ \al = 0, 1$), two vertical paths 
$ t \mapsto (\al, t)$, and all increasing paths from $ (0, 0) $ to $ (1, 1) $ in the ordered square, as 
exemplified in the left picture above.

	There is no finite set of generators, but we shall see that the fundamental category has 
only five non-trivial arrows (and four identities at the flexible points).

	The space is rigid. The generated d-space is the ordered square $ \uI^2$, whose d-paths 
are the increasing maps $ \bbI \to \bbI^2$; one of them is drawn in the second picture above.

	Similarly, in the {\em standard c-cube} $ \cI^n $ a path is controlled if and only if it is 
increasing and each of its projections covers $ [0, 1]$, or is constant at 0 or 1. Again, 
$ (\cI^n)\hat{\;} = \uI^n$.

\begin{small}

	For a product, the structure $ (X \ti Y)\hat{\;} $ is always finer than $\hat{X} \ti \hat{Y}$, and 
can be strictly finer. For instance, take $ X = \cI $ and $ Y = \cJ$, so that 
$ (\cI)\hat{\;} = (\cJ)\hat{\;} = \uI$; then $ (\cI \ti \cJ)\hat{\;} $ is finer than $ \uI^2$, and the diagonal 
$ \bbI \to \bbI^2 $ is controlled in $ \uI^2 $ but not in $ (\cI \ti \cJ)\hat{\;}$.

\end{small}

\Ndt (c) The {\em hybrid square} $ \cI \ti \uI $ will be important in the construction of the 
fundamental category. Here a path $ \bbI \to \bbI^2 $ is controlled if and only if it is increasing 
and its first projection is either surjective or constant at  0  or  1
%
% FIGURE 2.7.2, Hybrid controlled square
    \begin{equation} 
\xy <.5mm, 0mm>:
%horizontal POSITIONS of two squares: [0,30], [80,110]
%
% dummy pt above/left, below
(-20,22) *{}; (0,-22) *{};
% labels
(48,-10) *{\cI \ti \uI}; (125,-15) *{\sst{s}}; (119,0) *{\sst{t}};
% dashed square
@i@={(0,-15), (30,-15), (30,15), (0,15)},
s0="prev"  @@{;"prev";**@{--}="prev"};
% curved path and its arrow
(0, -10); (30, 10) **\crv{(5,5)&(20,0)},
\POS(15,2.5) \are+(4,1.2), 
% vertical paths
\POS(0,0) \ar+(0,10), \POS(30,-10) \ar+(0,10),
% axes of the plane
\POS(110,-10) \ar+(20,0),  \POS(115,-15) \ar+(0,20),
\endxy
    \label{2.7.2} \end{equation}

\begin{small}

	(All the horizontal paths $ s \mapsto (s, t_0) $ are controlled, but already belong to the 
family of increasing paths whose first projection is surjective.)

\end{small}

\medskip	The generated d-space is again $ \uI^2$; the flexible part $ D\{0, 1\} \ti \uI $ only allows 
increasing paths in the vertical edges.

\Ndt (d) The following example shows a sharp distinction between a c-structure $ X $ of the 
square and the generated d-space
%
% FIGURE 2.7.3, A controlled square and the generated d-space
    \begin{equation} 
\xy <.5mm, 0mm>:
%horizontal POSITIONS of two squares: [0,30], [90,120]
%
% dummy pt above/left, below
(-10,22) *{}; (0,-22) *{};
% labels, bullets
(45,-7) *{X}; (140,-7) *{\hat{X}}; (85,-15) *{\sst{p'}}; (126,-15) *{\sst{p''}};
(0,-15) *{\bu}; (30,-15) *{\bu}; (30,15) *{\bu};  (0,15) *{\bu}; 
% two squares
@i@={(0,-15), (30,-15), (30,15), (0,15)},
s0="prev"  @@{;"prev";**@{--}="prev"};
@i@={(90,-15), (120,-15), (120,15), (90,15)},
s0="prev"  @@{;"prev";**@{--}="prev"};
% paths and arrows
\POS(0,-15) \arl+(30, 30), \POS(0,15) \arl+(30, -30), 
\POS(90,-15) \arl+(15, 15), \POS(105,0) \arl+(15, -15), 
\POS(22,7) \are+(2,2), \POS(22,-7) \are+(2,-2), \POS(112,-7) \ar+(2,-2), 
\endxy
    \label{2.7.3} \end{equation}

	The c-paths of $ X $ are generated by two diagonal paths, $ t \mapsto (t, t) $ and 
$ t \mapsto (t, 1 - t)$; the flexible points are the four vertices of the square, and it is easy to 
guess that the fundamental category $ \uPi_1(X) $ will only have two non-trivial arrows, the 
marked ones. But the fundamental category $ \uPi_1(\hat{X})) $ of the generated d-space 
has many others, including two new arrows linking the vertices: for instance, the arrow 
$ p' \to p'' $ displayed in the right figure above.
	
	Clearly, one cannot model a crossing of railways by a d-space. Within c-spaces, one
 can make $ X $ reversible adding the opposite c-paths.

%%%
%%%
\section{Critical processes and models}\label{s3}
	We now start from `critical processes', trying to model them by c-spaces built with the 
previous ones, by limits and colimits.

\subsection{On-off controllers and elementary hysteresis}\label{3.1}
We consider a system meant to regulate a certain variable, either opposing its rising, or 
helping it, or working both ways. An elementary hysteresis process, or `hysteron', generally 
behaves in the first way -- counteracting the effect.

\Ndt  (a) {\em Reacting controller.} We begin considering a cooling device counteracting the rising 
of temperature, with a thermostat set at temperature $T_0$ and a tolerance interval $[T_1, T_2]$.

	In the following picture the horizontal axis represents the temperature, and the vertical 
axis denotes two states: at level 0 the system is off, at level 1 it is on
%
% FIGURE 3.1.1  On-off controller, from Fig. 0.2.1
    
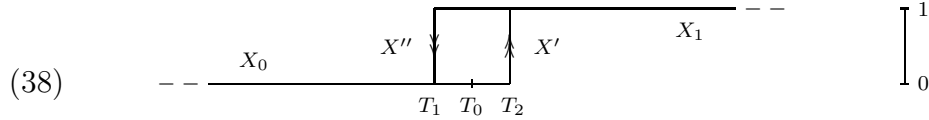
\begin{equation}
\xy <.5mm, 0mm>:
% labels + dummy pts left/above, below
(12,6) *{\sst{X_0}}; (50,10) *{\sst{X''}}; (90,10) *{\sst{X'}};  (128,14) *{\sst{X_1}}; 
(59,-6) *{\sst{T_1}}; (70,-6) *{\sst{T_0}}; (81,-6) *{\sst{T_2}}; 
(190,0) *{\sst{0}}; (190,20) *{\sst{1}};  
(-30,30) *{}; (-30,-20) *{};
% branches X_0, X_1
\POS(0,0) \arl+(80,0),  \POS(60,20) \arl+(80,0),
\POS(-13,0) \arld+(10,0),  \POS(143,20) \arld+(10,0),  
\POS(70,-1) \arl+(0,2),
% branches X', X'' and double arrows
\POS(60,0) \arl+(0,20),  \POS(80,0) \arl+(0,20),
\POS(60,12) \are+(0,-5),  \POS(80,8) \are+(0,5),
%vertical axis
\POS(185,0) \arl+(0,20), \POS(184,0) \arl+(2,0), \POS(184,20) \arl+(2,0), 
\endxy
    \label{3.1.1} \end{equation}

	On the left branch $X_0$ the system is off; if the temperature grows to $T_2 $ the device jumps 
to state 1; then, if the temperature cools to $T_1$, it goes back to state 0.

	The support $ |X| $ of our model is a one-dimensional subspace of $ \bbR^2$, the union 
of the supports of the following c-spaces
    \begin{equation*} \begin{array}{ccc}
X_0  =  [0, T_2] \times \{0\},   &\q&   X_1  =  [T_1, + \infty[ \, \times \{1\},
\\[2pt]
X'  =  \{T_2\} \times \, \cI,   &&   X''  =  \{T_1\} \times \cI\op.
    \label{} \end{array} \end{equation*}

	The c-structure of $ X $ is generated by the c-structures of:
	
\Ndt - $ X_0, X_1$, natural intervals, where the temperature can vary, 

\Ndt - $ X', X''$, one-jump c-intervals, where the state of the system varies.

\begin{small}

\smallskip	 One could also use the plane with the terminal c-structure $ \rc_f\bbR^2 $ produced by 
the topological embedding $ f\c X \to \bbR^2$; the c-paths are those of $ X$.

\smallskip	 A hysteretic process is generally studied as a functional operator that turns a piecewise 
monotone input function into an output function (of time, in both cases): see \cite{BrS}, Chapter 2. 
Here the input is the temperature function, while the output values are the 
states 0, 1. This analysis presents some indetermination and failure of continuity at the critical 
temperatures $T_1, T_2$, as discussed in \cite{BrS}, Example 2.1.1.
	
\end{small}

\Ndt (b) {\em Cooperating controller.} A heating system supports the rising of temperature. It can 
be modelled by the opposite c-space: in state 0 the system is on; see the lower half of the following diagram.

\Ndt (c) {\em Dual controller.} Combining both models we can represent a heating and cooling 
system, like a heat pump. The system is meant to keep the temperature in an interval 
$[T_1, T_2]$, with a lower tolerance $[T_1, T'_1]$ and an upper tolerance $[T'_2, T_2]$ 
(disjoint intervals); the vertical axis denotes now three states: at level 0 the system is off, at 
level -1 the heating is on, at level 1 the cooling device is on
%
% FIGURE 3.1.2  Dual on-off controller
    
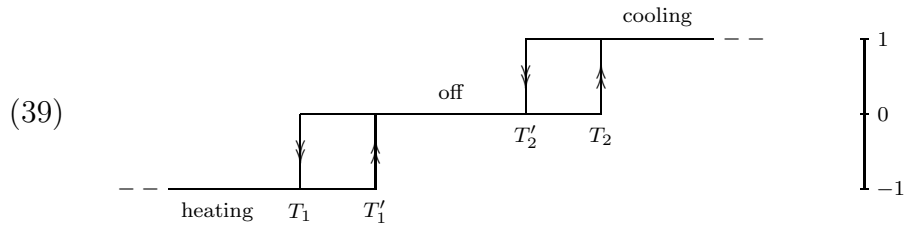
\begin{equation}
\xy <.5mm, 0mm>:
% labels + dummy pts left/above, below
(-22,-26) *{\sst{\text{heating}}}; (40,6) *{\sst{\text{off}}}; (95,26) *{\sst{\text{cooling}}}; 
(0,-26) *{\sst{T_1}}; (20,-26) *{\sst{T'_1}}; 
(60,-6) *{\sst{T'_2}}; (80,-6) *{\sst{T_2}}; 
(157,-20) *{\sst{-1}}; (155,0) *{\sst{0}}; (155,20) *{\sst{1}};  
(-30,30) *{}; (-30,-40) *{};
% horizontal lines
\POS(-35,-20) \arl+(55,0),  \POS(0,0) \arl+(80,0),  \POS(60,20) \arl+(50,0),
\POS(-48,-20) \arld+(10,0),  \POS(113,20) \arld+(10,0),  
% vertical lines and double arrows
\POS(0,-20) \arl+(0,20),  \POS(20,-20) \arl+(0,20),
\POS(60,0) \arl+(0,20),  \POS(80,0) \arl+(0,20),
\POS(0,-8) \are+(0,-5),  \POS(20,-12) \are+(0,5),
\POS(60,12) \are+(0,-5),  \POS(80,8) \are+(0,5),
%vertical axis
\POS(150,-20) \arl+(0,40), 
\POS(149,-20) \arl+(2,0) \POS(149,0) \arl+(2,0), \POS(149,20) \arl+(2,0), 
\endxy
    \label{3.1.2} \end{equation}

\subsection{Controlling two variables}\label{3.2}
We deal now with a pair of on-off controllers, acting on two independent variables

\Ndt (a) {\em Two reacting controllers}. We start from two copies $ X, Y $ of the c-space drawn 
in \eqref{3.1.1}
%
% FIGURE 3.2.1  Two on-off controllers
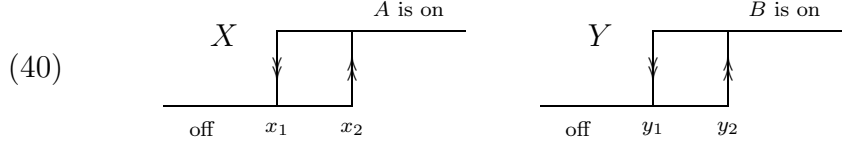
    \begin{equation}
\xy <.5mm, 0mm>:
% labels + dummy pts left/above, below
(-14,9) *{X}; (86,9) *{Y}; 
(-20,-16) *{\sst{\text{off}}}; (35,16) *{\sst{A \text{ is on}}}; 
(80,-16) *{\sst{\text{off}}}; (135,16) *{\sst{B \text{ is on}}}; 
(0,-16) *{\sst{x_1}}; (20,-16) *{\sst{x_2}}; 
(100,-16) *{\sst{y_1}}; (120,-16) *{\sst{y_2}}; 
(-40,20) *{}; (-30,-30) *{};
% horizontal lines
\POS(-30,-10) \arl+(50,0), \POS(0,10) \arl+(50,0), 
\POS(70,-10) \arl+(50,0), \POS(100,10) \arl+(50,0),
% vertical lines and double arrows
\POS(0,-10) \arl+(0,20),  \POS(20,-10) \arl+(0,20),
\POS(100,-10) \arl+(0,20),  \POS(120,-10) \arl+(0,20),
\POS(0,2) \are+(0,-5),  \POS(20,-2) \are+(0,5),
\POS(100,2) \are+(0,-5),  \POS(120,-2) \are+(0,5),
\endxy
    \label{3.2.1} \end{equation}

	In $ X $ a device $ A $ controls the variable $ x$, countering its rising; in $ Y $ the device 
$ B $ acts similarly on the variable $ y$.
	
	The cartesian product $ X \ti Y $ models the combined system. Its support is a subspace of 
$ \bbR^4$, but we draw it in $ \bbR^3$, with four states on the vertical axis: at 0 both devices 
are off, at 1 only $ A $ is on, at 2 only $ B $ is on, at 3 both $ A $ and $ B $ are on
%
% FIGURE 3.2.2  Two combined on-off controllers
    
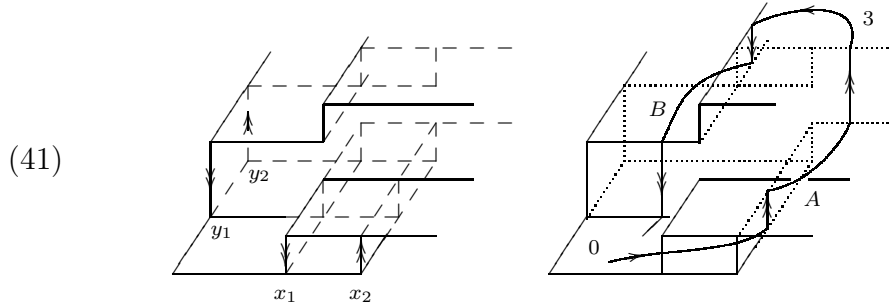
\begin{equation}
\xy <.5mm, 0mm>:
% LEFT FIGURE
% labels + dummy pts left/above, below
(30,-36) *{\sst{x_1}}; (50,-36) *{\sst{x_2}}; (13,-19) *{\sst{y_1}}; (23,-4) *{\sst{y_2}}; 
(-10,50) *{}; (0,-40) *{};
% front horizontal lines, with vertical jumps and arrows
\POS(0,-30) \arl+(50,0), \POS(30,-20) \arl+(40,0), 
\POS(30,-30) \arl+(0,10),  \POS(50,-30) \arl+(0,10),
\POS(30,-23) \are+(0,-5),  \POS(50,-26) \are+(0,5),
% left slanting lines, with vertical jumps and arrows
\POS(0,-30) \arl+(10,15), \POS(10,-15) \arld+(10,15), \POS(10,5) \arl+(16,24), 
\POS(10,-15) \arl+(0,20),  \POS(20,0) \arld+(0,20),
\POS(10,-3) \are+(0,-5),  \POS(20,8) \are+(0,5),
% middle horizontal lines, with vertical jumps
\POS(10,-15) \arl+(20,0), \POS(30,-15) \arld+(30,0), \POS(40,-5) \arl+(40,0), 
\POS(10,5) \arl+(30,0), \POS(40,15) \arl+(40,0), 
\POS(20,0) \arld+(50,0), \POS(50,10) \arld+(40,0), 
\POS(20,20) \arld+(50,0), \POS(50,30) \arld+(40,0), 
\POS(70,0) \arld+(0,10), \POS(50,20) \arld+(0,10), \POS(70,20) \arld+(0,10), 
\POS(40,-15) \arld+(0,10), \POS(60,-15) \arld+(0,10),
\POS(40,5) \arl+(0,10), \POS(50,0) \arld+(0,10), 
% middle slanting lines
\POS(30,-30) \arld+(20,30), \POS(50,-30) \arl+(8,12), \POS(55,-22.5) \arld+(14,21), 
\POS(30,-20) \arl+(10,15), \POS(40,-5) \arld+(8,12), \POS(50,-20) \arld+(20,30),
\POS(40,15) \arl+(12,18), \POS(40,5) \arld+(12,18), 
% RIGHT FIGURE: horizontal + 100
% labels
\POS(0,0),
(112,-23) *{\sst{0}}; (170,-10) *{\sst{A}}; (185,38) *{\sst{3}}; (129,14) *{\sst{B}}; 
% front horizontal lines, with vertical jumps
\POS(100,-30) \arl+(50,0), \POS(130,-20) \arl+(40,0),
\POS(130,-30) \arl+(0,10),  \POS(150,-30) \arl+(0,10),
% left slanting lines, with vertical jumps
\POS(100,-30) \arl+(10,15), \POS(110,-15) \arlp+(10,15), \POS(110,5) \arl+(16,24), 
\POS(120,0) \arlp+(0,20), \POS(110,-15) \arl+(0,20),
% middle horizontal lines, with vertical jumps
\POS(110,-15) \arl+(22,0), \POS(140,-5) \arl+(24,0), \POS(169,-5) \arl+(11,0),    %%% 
\POS(110,5) \arl+(30,0), \POS(140,15) \arl+(20,0), 
\POS(120,0) \arlp+(50,0), \POS(170,0) \arlp+(0,10),  \POS(170,10) \arlp+(20,0),
\POS(120,20) \arlp+(50,0), \POS(150,20) \arlp+(0,10),  \POS(170,20) \arlp+(0,10), 
\POS(150,30) \arlp+(40,0), \POS(140,5) \arl+(0,10),
% middle slanting lines
\POS(150,-30) \arl+(6.4,9.6), \POS(155,-22.5) \arlp+(14,21),
\POS(130,-20) \arl+(10,15), \POS(150,-20) \arlp+(20,30),
\POS(140,15) \arl+(18,27), \POS(140,5) \arlp+(18,27),
%the path, jumps and arrows
\POS(158, -18) \arl+(0,10), \POS(158, -13) \are+(0,4), 
\POS(180, 10) \arl+(0,20), \POS(180, 20) \are+(0,4), 
\POS(154, 36) \arl+(0,-10), \POS(154, 31) \are+(0,-4), 
\POS(130, 5) \arl+(0,-20), \POS(130, -5) \are+(0,-4), 
\POS(130, -15) \arl+(-5,-5),
\POS(121, -26) \ar+(4,.8), \POS(170, 40) \ar+(-3,0), 
\POS(0,0),
%the path, curved parts
(116, -27); (158, -18) **\crv{(130,-23)&(150,-25)},
(158, -8); (180, 10) **\crv{(170,-6)&(178,5)},
(180, 30); (154, 36) **\crv{(185,45)&(160,40)},
(154, 26); (130, 5) **\crv{(135,22)&(135,15)},
\endxy
    \label{3.2.2} \end{equation}

	The new c-space can be obtained as follows. We put on $ \bbR^3 $ the terminal c-structure for the topological map
    \begin{equation}
f\c X \ti Y \to \bbR^3,   \qq   f(x, s, y, t)  =  (x, y, s + 2t),
    \label{3.2.3} \end{equation}
and we use this c-space $ \rc_f\bbR^3$. Equivalently, we can use the c-subspace $ Z = \Im f $ of 
the previous structure: they have the same c-paths. $ Z $ is contained in four parallel planes, at 
level 0, 1, 2, 3.

	The right figure above represents a c-path in $ Z$. It starts at level 0, with both variables 
below the active thresholds $ x_2, y_2$, and both increasing; when the variable $ x $ reaches 
$ x_2$, system $ A $ jumps on, counteracting it. Both variables are still growing; when $ y $ 
attains $ y_2$, device $ B $ also activates and the process is in state 3. Then the variable $ x $ 
decreases below $ x_1 $ and $ A $ goes off, while $ B $ is still on, in state 2. Finally 
also the variable $ y $ is brought below $ y_1 $ and both devices are off, at level 0.

\Ndt (b) The opposite case of two cooperating controllers is modelled by the opposite c-space 
$ Z\op$. The mixed case, with $ A $ cooperating with variable $ x $ and $ B $ counteracting 
variable $ y$, is also of interest.

	An air-supported dome can give examples of both cases. In winter time, the compressor 
would rise pressure and the heating would rise temperature; in summer time, the compressor 
would work the same way while air-conditioning would reduce temperature.

\subsection{The threshold effect and siphon structures. }\label{3.3}
In the threshold effect a variable $v$ can vary in the interval $[v_0, v_1]$; reaching the 
{\em threshold} $v_1$ it jumps down to its least value $v_0$. Various processes of this type 
are listed in the Introduction, Subsection \ref{0.4}. Here we consider two structures on the 
interval $\bbI$ that can model such a process.

\Ndt (a) {\em The growing siphon}. We denote as ${\rm c}_S\bbI$ the standard interval with 
the c-structure generated by all the increasing maps $\bbI \to \bbI$ and the reversion $r(t) = 1 - t$
%
% FIGURE 3.3.1, the growing siphon
    \begin{equation}
\xy <.5mm, 0mm>:
% labels + dummy pt above, below
(0,-6) *{\sst{0}}; (50,-6) *{\sst{1}}; 
(0,10) *{}; (0,-20) *{};
% the interval
\POS(0,0) \arl+(50,0), \POS(0,-1) \arl+(0,2), \POS(50,-1) \arl+(0,2), 
\POS(0,0), (3,-1); (47,-1) **\crv{~*=<5pt>{.}(10,-6)&(40,-6)};     % dotted curve 
\POS(23,0) \ar+(4,0), \POS(25.5,-4.7) \are+(-4,0), 
\endxy
    \label{3.3.1} \end{equation}

	A controlled path can only increase between 0 and 1; reaching 1, either it stays there or 
jumps down to 0. Point 1 is future critical, point 0 is past critical, and there are no (bilateral) 
critical points. The generated d-structure $({\rm c}_S\bbI)\hat{\;}$ has two generators, 
$ \id \bbI $ and $r$: it is the reversible d-interval $ \bbI^\sim$ of \ref{2.4}(c); the flexible part is 
$ \uI$.

	For concreteness, we refer to a hydraulic system consisting of a water basin filled by a 
source; water can only get out by a siphon tube, as in the figure below, so that its level $ v $ in 
the basin will grow up to the upper part of the tube, marked 1, and then flow out until the level 
reaches the lower opening of the tube in the basin, marked 0 (the diameter of the tube is
overlooked)
%
% FIGURE 3.3.2  a hydraulic system
    
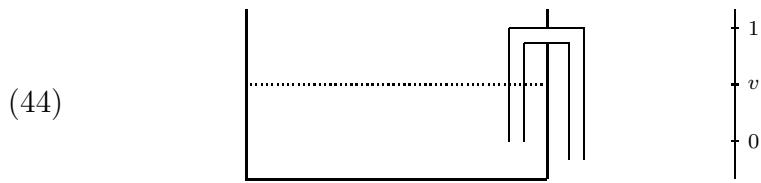
\begin{equation}
\xy <.5mm, 0mm>:
% labels + dummy pts left/above, below
(155,-10) *{\sst{0}}; (155,5) *{\sst{v}}; (155,20) *{\sst{1}};  
(10,35) *{}; (10,-25) *{};
% basin (double lines), water level
\POS(20,-20) \arl+(0,45), \POS(20,-20) \arl+(80,0), \POS(100,-20) \arl+(0,36), 
\POS(100,20) \arl+(0,5), 
\POS(20.3,-20) \arl+(0,45), \POS(20,-20.3) \arl+(80,0), \POS(100.3,-20) \arl+(0,36), 
\POS(100.3,20) \arl+(0,5), 
\POS(20,5) \arlp+(80,0),
% siphon tube
\POS(90,-10) \arl+(0,30), \POS(90,20) \arl+(20,0), \POS(110,20) \arl+(0,-35), 
\POS(94,-10) \arl+(0,26), \POS(94,16) \arl+(12,0), \POS(106,16) \arl+(0,-31), 
%vertical axis
\POS(150,-20) \arl+(0,45), 
\POS(149,-10) \arl+(2,0) \POS(149,5) \arl+(2,0), \POS(149,20) \arl+(2,0), 
\endxy
    \label{3.3.2} \end{equation}

\Ndt (b) {\em The oscillating siphon}. A more complex model ${\rm c}'_S\bbI $ allows the c-paths to 
decrease in the semiopen interval $ [0, 1[$. There are three kinds of generators of the c-paths:

\Ndt - the increasing maps $ \bbI \to \bbI$,

\ndt - the decreasing maps $ a\c \bbI \to \bbI $ with image in $ [0, 1[$, i.e. $ a(0) < 1$,

\ndt - the reversion $ r(t) = 1 - t$.

\skp	Here a controlled path is piecewise monotone; whenever it reaches 1, either it stays there or 
it decreases to 0 (and possibly goes on). Point 1 is still future critical, but 0 is no more past critical. 
The generated d-structure is again the reversible d-interval $ \bbI^\sim$, the flexible part 
allows the piecewise monotone paths which can only reach 1 as the terminal endpoint.

	In the hydraulic system previously described, the basin can now let out water by other 
openings or evaporation. Other processes considered in \ref{0.4} are also better fitted by this 
model; for instance the electric potential at a neural membrane can increase and decrease; 
reaching the threshold value the impulse is emitted.

\subsection{Transport networks and labelled graphs}\label{3.4}
Transport networks are usually modelled within graph theory. They can also be modelled by 
c-spaces, as we have already seen in various examples; this would allow to combine them 
with planar or three-dimensional regions.

\Ndt (a) The following is a model of a road that contains a dual-carriage section (for left-hand 
drive one would turn the picture upside-down)
%
% FIGURE 3.4.1, a road with a dual-carriage section
    
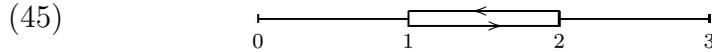
\begin{equation}
\xy <.5mm, 0mm>:
% labels + dummy pt above, below
(0,-6) *{\sst{0}}; (40,-6) *{\sst{1}}; (80,-6) *{\sst{2}}; (120,-6) *{\sst{3}};
(0,15) *{}; (0,-25) *{};
% the line
\POS(0,0) \arl+(40,0), \POS(40,-2) \arl+(40,0), \POS(40,2) \arl+(40,0), 
\POS(80,0) \arl+(40,0), 
\POS(0,-1) \arl+(0,2), \POS(40,-2) \arl+(0,4), 
\POS(80,-2) \arl+(0,4), \POS(120,-1) \arl+(0,2),
\POS(61,-2) \ar+(4,0), \POS(61,2) \ar+(-4,0), 
\endxy
    \label{3.4.1} \end{equation}

	The d-space $X$ is the quotient $Y/R$ of the sum $Y = X_1 + ... + X_4$ of the 
following d-spaces
    \begin{equation*} \begin{array}{ccc}
X_1  =  [0, 1],    &\q&    X_2  =  \uw[1, 2],
\\[5pt]
X_3  =  \uw[1, 2]\op,    &&    X_4  =  [2, 3],
    \label{} \end{array} \end{equation*}
modulo the equivalence relation that identifies the three points 1 (of $X_1, X_2, X_3$) and 
(separately) the three points 2 (of $ X_2, X_3, X_4$). A path in $ X $ is directed if and only if it is 
a general concatenation of projections of d-paths in the various $ X_i$. For instance, to go 
from 0 to the point $ (1/2)_3 $ of $ X_3 $ we must (at least) reach $ 2 $ along $ X_1 $ and 
$X_2 $ and then come back along $ X_3$; there are infinitely many longer paths.
 
\Ndt (b) Similarly, one can construct a one-dimensional c-space $X$ as a realisation of a 
 {\em labelled graph}, in the sense of a multigraph whose edges are labelled with additional 
 information on direction and critical properties. We have already drawn many linear examples, 
 and some planar ones in \ref{3.1}.

	As above, the c-space $X$ can be obtained as a quotient of a sum of intervals with the 
appropriate c-structure: natural intervals $\bbI$ when unlabelled, standard d-intervals
$\uI$ when labelled by a single arrow, standard c-intervals $\cI$ when labelled by a double 
arrow, etc.

	This can represent a transport networks, where some routes are one-way and others are 
no-stop -- as in an underground section, or a motorway tunnel, or an airline route. 
The model can be further enriched, adding delays (for a stop sign), etc. Or higher dimensional 
regions, as we were suggesting.

\subsection{Point-like variations}\label{3.5}
In the examples of this section one can often form a `slightly' different model, using a general 
procedure: if $X$ is a c-space and $A \sub |X|$, one builds a finer c-space on $|X|$ excluding all 
the previous c-paths that have an endpoint in $A$.

	Thus, in the model of the heat controller described in \ref{3.1}(a), one can omit 
the paths that start or end at $(0, T_2)$ or $(1, T_1)$. Similarly, in the siphon-interval 
\ref{3.3}(a) one can rule out the paths starting or ending at 1. In both cases we are forcing the 
jump at these points, which become critical and non-flexible. This can be appreciated, but the 
new models are more complicated and their fundamental category will also be.

	In our opinion the choice between such variations is merely a theoretical issue, 
that should be based on the results one can obtain. 
In the same way as, if we model a thin rod by the interval $[0, 1]$, it is Mathematics rather than 
experience that leads us to use an interval of the real line, instead of the rational line where 
the classical results on continuous functions and differential equations would fail.

%%%
%%%

\end{document}